\documentclass[11pt,reqno]{amsart}
\usepackage{graphicx} 
\usepackage{hyperref}
\usepackage{fourier}
\usepackage{fullpage}
\usepackage{amsmath,amssymb,amsthm}
\usepackage[shortlabels]{enumitem}
\usepackage[mathscr]{euscript}
\usepackage{dutchcal}
\usepackage{upgreek}
\usepackage{comment}
\usepackage{csquotes}
\usepackage{tikz}
\allowdisplaybreaks
\providecommand{\noopsort}[1]{}

\makeatletter
\def\resetMathstrut@{%
  \setbox\z@\hbox{%
    \mathchardef\@tempa\mathcode`\(\relax
    \def\@tempb##1"##2##3{\the\textfont"##3\char"}%
    \expandafter\@tempb\meaning\@tempa \relax
  }%
  \ht\Mathstrutbox@1.2\ht\z@ \dp\Mathstrutbox@1.2\dp\z@
}
\makeatother

\newtheorem{theorem}{Theorem}
\newtheorem{lemma}[theorem]{Lemma}

\newtheorem{remark}[theorem]{Remark}

\newcommand{\R}{\mathbb{R}}

\newcommand{\N}{\mathbb{N}}

\renewcommand{\le}{\leqslant}
\renewcommand{\ge}{\geqslant}
\renewcommand{\leq}{\leqslant}

\renewcommand{\setminus}{\smallsetminus}
\renewcommand{\subset}{\subseteq}

\newcommand{\eqdef}{\stackrel{\mathrm{def}}{=}}

\newcommand{\X}{\mathbf X}
\newcommand{\Y}{\mathbf Y}

\newcommand{\bZ}{\mathbf{Z}}

\newcommand{\ud}[0]{\,\mathrm{d}}

\newcommand{\sfA}{\mathsf{A}}

\newcommand{\n}{\{1,\ldots,n\}}
\newcommand{\sfC}{\mathsf{C}}

\newcommand{\GL}{\mathsf{GL}}
\newcommand{\BM}{\mathrm{BM}}

\renewcommand{\nu}{\upnu}

\newcommand{\og}{\othergamma}
\newcommand{\sfT}{\mathsf{T}}
\newcommand{\sfS}{\mathsf{S}}
\newcommand{\sfM}{\mathsf{M}}

\newcommand{\conv}{\mathrm{conv}}
\newcommand{\Id}{\mathsf{I}}
\newcommand{\Tr}{\mathrm{trace}}

\newcommand{\sR}{\mathscr{R}}

\usepackage{xcolor}
\definecolor{maroon}{HTML}{AF3235}

\title{On the Banach--Mazur radius of $\ell_\infty^n$}

\author{Omer Friedland}
\address{Omer Friedland. Institut de Math\'ematiques de Jussieu, Sorbonne Universit\'e, 4 place Jussieu, 75005 Paris, France.}
\email{omer.friedland@imj-prg.fr}

\author{Assaf Naor}
\address{Assaf Naor. Mathematics Department\\ Princeton University\\ Fine Hall, Washington Road, Princeton, NJ 08544-1000, USA.}
\email{naor@math.princeton.edu}

\author{Pierre Youssef}
\address{Pierre Youssef. Center for Interdisciplinary Data Science and AI, NYUAD Research Institute, UAE \& Division of Science, NYU Abu Dhabi, Abu Dhabi, UAE \& Courant Institute of Mathematical Sciences, New York University, New York, USA.}
\email{yp27@nyu.edu}

\thanks{A.~N.~is supported by NSF grant DMS-2453936 and a Simons Investigator award. P.~Y.~is supported by the NYUAD Center for Interdisciplinary Data Science
\& AI (CIDSAI), funded by Tamkeen under the NYUAD Research Institute Award CG016.} 

\date{}
\begin{document}

\maketitle
\vspace{-2em}
\begin{abstract} 
We prove that any $n$-dimensional normed space is at Banach--Mazur distance $O(n^{2/3})$ from $\ell_\infty^n$. 
\end{abstract}

\section{Introduction}

Fix $n\in \N$ and an origin-symmetric convex body $K\subset \R^n$. Thanks to John's theorem~\cite{MR30135} there exists an origin-symmetric ellipsoid $\mathcal{E}\subset \R^n$ and  $0<R\le \sqrt{n}$ such that $\mathcal{E}\subset K\subset R\mathcal{E}$. E.g.~$K=[-1,1]^n$ shows that this is optimal.  Understanding the analogous  phenomenon in which one wishes to approximate any  origin-symmetric convex body by  a parallelepiped  is a major  mystery. Specifically, let $\sR_{\BM}(\ell_\infty^n)$ be the smallest $R>0$ such that every origin-symmetric convex body $K\subset \R^n$ admits an origin-symmetric  parallelepiped  $\mathcal{P}\subset \R^n$ with  $\mathcal{P}\subset K\subset R\mathcal{P}$.  Gr\"unbaum~\cite[Remark~6(ii)]{Gru60} was the first to publish the natural question of evaluating $\sR_{\BM}(\ell_\infty^n)$.  Pe{\l}czy{\'n}ski conjectured~\cite{Pelczynski1984} that  $\lim_{n\to\infty} \sR_{\BM}(\ell_\infty^n)/\sqrt{n}=\infty$.

In~\cite{MR1081810}, Szarek confirmed the  above conjecture of Pe{\l}czy{\'n}ski by proving that $\sR_{\BM}(\ell_\infty^n)\gtrsim \sqrt{n} \log n$.\footnote{In  addition to the usual $O(\cdot)$ notation, we use the following standard asymptotic notation:  $A\lesssim B$ and $B\gtrsim A$ stand for $A\le c B$ with $c$ a universal constant,  $A\asymp B$ stands for
$A\lesssim B\lesssim A$. } Tikhomirov~\cite{MR3900705} and the first named author~\cite{Fri26} improved this to  $\sR_{\BM}(\ell_\infty^n)\ge n^{5/9}/(\log n)^{O(1)}$ and  $\sR_{\BM}(\ell_\infty^n)\ge n^{4/7}/(\log n)^{O(1)}$, respectively. The most recent progress on  lower bounds for  $\sR_{\BM}(\ell_\infty^n)$ is due to~\cite{Fri26-5/6} and independently Hmadi~\cite{Hma26}, where  $\sR_{\BM}(\ell_\infty^n)\gtrsim n^{5/8}/(\log n)^{O(1)}$ is obtained. 

Auerbach’s lemma~\cite[page~238]{Ban32} gives $\sR_{\BM}(\ell_\infty^n)\le n$. The first asymptotic improvement over this classical upper bound was by Bourgain and Szarek~\cite{MR947820}, who proved that  $\sR_{\BM}(\ell_\infty^n)\le n\exp({-c\sqrt{\log n}})$ for some universal constant $c>0$. Szarek and Talagrand~\cite{MR1008718} and Giannopoulos~\cite{MR1353450} improved this to $\sR_{\BM}(\ell_\infty^n)\lesssim  n^{7/8}$ and $\sR_{\BM}(\ell_\infty^n)\lesssim  n^{5/6}$, respectively (see also~\cite{MR3164527} for a simplified approach and a smaller bound on the implicit universal constant factor). The present work is devoted to proving:

\begin{theorem}\label{thm:main} $\mathscr{R}_{\BM}(\ell_\infty^n)\lesssim n^{\frac23}$ for every $n\in \N$. 
\end{theorem}

The above topic is naturally expressed and studied in the context of the  Banach--Mazur compactum. Given isomorphic  Banach spaces $\X,\Y$,   the infimum of $\|\sfT\|_{\X\to\Y}\|\sfT^{-1}\|_{\Y\to\X}$ over all isomorphisms $\sfT:\X\to \Y$ is called their Banach--Mazur distance, and it is denoted $d_{\BM}(\X,\Y)$.  The Banach--Mazur radius $\sR_{\BM}(\bZ)$ of a Banach space $\bZ$ is the supremum of $d_{\BM}(\bZ,\X)$ over all  Banach spaces $\X$ that are isomorphic to $\bZ$. For $\bZ=\ell_\infty^n$, this notation coincides with the above defined parallelepiped squeezing radius $\sR_{\BM}(\ell_\infty^n)$, as seen by taking $K=B_\X=\{x\in \X:\ \|x\|_\X\le 1\}$ to be the unit ball of a normed space $\X=(\R^n,\|\cdot\|_\X)$, and noting that a parallelepiped is the same as $\sfT B_{\ell_\infty^n}=\sfT[-1,1]^n$ for some $\sfT\in \GL_n(\R)$.  Observe also  the duality  $\sR_{\BM}(\bZ)=\sR_{\BM}(\bZ^*)$, whence we will use below $\sR_\BM(\ell_\infty^n)=\sR_{\BM}(\ell_1^n)$.   

Despite intensive investigations of Banach--Mazur distances over many decades and a lot of accumulated knowledge (see e.g.~the monograph~\cite{MR995162}), the Banach--Mazur radius of all but a few spaces has not been determined (even up to lower order factors);   $\sR_\BM(\ell_\infty^n)$ is the most prominent such example.\footnote{While the present work concerns  finite dimensional normed spaces, it is worthwhile to recall the following longstanding question about Banach--Mazur radii in infinite dimensions. One is tempted  to expect that $\sR_{\BM}(\X)=\infty$ for every infinite dimensional Banach space $\X$, but this turns out to be far from obvious; the question whether this holds was posed by Sch\"affer~\cite{Sch76}.  Johnson and Odell proved~\cite{JO05} that indeed $\sR_{\BM}(\X)=\infty$  for every infinite dimensional separable Banach space $\X$. Whether or not the same holds for every nonseparable Banach space remains open; see also~\cite{God10,Joh18}.} The resolved instances   are:  $\sR_{\BM}(\ell_2^n)=\sqrt{n}$ by John's theorem~\cite{MR30135}; Gluskin's theorem~\cite{Glu81} furnishes for every $1\le \alpha\le \sqrt{n}$ a distribution over random normed spaces $\mathbf{G}_\alpha=(\R^n,\|\cdot\|_{\mathbf{G}_\alpha})$ for which $\sR_{\BM}(\mathbf{G}_\alpha)\asymp \alpha\sqrt{n}$ with high probability (Remark~\ref{rem:gluskin} below explains how this follows quickly from~\cite{Glu81}, where only the case $\alpha=\sqrt{n}$ is stated, famously showing that the diameter of the Banach--Mazur compactum over $\R^n$ is of order $n$);  Bourgain and Szarek proved~\cite{MR947820} that $\sR_{\BM}(\ell_2^n\oplus \ell_\infty^n)=  \sR_{\BM}(\ell_2^n\oplus \ell_1^n)\asymp \sqrt{n}$, thus demonstrating (and answering negatively a conjecture of Pe{\l}czy{\'n}ski~\cite{Pelczynski1984}) that there are non-Hilbertian asymptotic centers of the Banach--Mazur compactum.\footnote{Formally, this  list is not exhaustive as it can be modified in inessential ways: other Gluskin-like random spaces  can be used, and  the dimensions of the direct summands in the Bourgain--Szarek theorem can be any other universal proportions of $n$. } 

The obvious question that remains is to determine $\sR_{\BM}(\ell_\infty^n)$. A $2/3$ exponent  has an intrinsic meaning (in a precise sense that we do  not discuss  herein)  for the  lower bound example that the aforementioned works~\cite{MR1081810,MR3900705,Fri26,Fri26-5/6,Hma26} consider (and natural generalizations thereof), as well as for the approach of the present article. So, we have reasons  to suspect that Theorem~\ref{thm:main} might be  optimal up to lower order factors, but at this juncture it would be speculative to conjecture that this indeed holds. 

\subsection*{On the proof of Theorem~\ref{thm:main}} Here we will describe (and motivate) the main steps by which  Theorem~\ref{thm:main} is proved. Fix $n\in \N$. Given a normed space $\X=(\R^n,\|\cdot\|_\X)$, the goal is to demonstrate that $d_\BM(\ell_1^n,\X)\lesssim n^{2/3}$.  As the unit ball of $\ell_1^n$ equals $\conv\{\pm e_1,\ldots, \pm e_n\}$, where $e_1,\ldots,e_n$ is the standard basis of $\R^n$, this amounts to a basis selection problem in which one seeks a well-conditioned basis that is suitably adapted to $\X$. The previous approaches pioneered by Bourgain and Szarek~\cite{MR947820} chose almost orthogonal contact points of $B_\X$ and the ellipsoid of minimum volume that contains it, and then added to  them orthonormal vectors so as to obtain a basis of $\R^n$, and  this procedure was performed iteratively by Giannopoulos~\cite{MR1353450} (the proof in~\cite{MR3164527}  utilizes for this the more recent tools of~\cite{BSS12,SS12}, as surveyed in~\cite{Nao12}).

The precedent  and inspiration for the present article is~\cite{NY18}, which proposes a framework for bounding $d_{\BM}(\ell_1^n,\X)$ by optimizing directly over full coordinate systems; this is part of a series of works~\cite{NY17-published,NY18,NY19} by the last two named authors that are devoted to related issues. Here we bring this project to fruition with further ideas that are described below. Our reasoning  turns out to be  self-contained without any need to appeal to~\cite{NY18}; this is beneficial also because~\cite{NY18} is still a preprint (that has been circulated to experts over the years), which will appear elsewhere as it is devoted to other goals.

The following parameter is a close relative of  the one that~\cite{NY18} proposed for bounding $d_\BM(\ell_1^n,\X)$:
\begin{equation}\label{eq:def og}
\og_\X\eqdef \min_{\substack{\sfT\in \GL_n(\R)\\ \sfT B^n\supseteq B_\X}}\bigg(\frac{1}{n}\sum_{i=1}^n \|\sfT e_i\|_\X^2\bigg)^{\frac12}=\min_{\sfT\in \GL_n(\R)}\bigg(\frac{1}{n}\sum_{i=1}^n \|\sfT e_i\|_\X^2\bigg)^{\frac12}\|\sfT^{-1}\|_{\X\to \ell_2^n},
\end{equation}
where $B^n=B_{\ell_2^n}$. The minimization in~\eqref{eq:def og} selects at the outset the basis $\sfT e_1,\ldots,\sfT e_n$, rather than obtaining a basis  in successive stages. The relevance of $\og_\X$ to Theorem~\ref{thm:main} stems from the following simple estimate:
\begin{equation}\label{eq:og bouds BM}
d_{\BM}(\ell_1^n,\X)\le \og_\X\sqrt{n}.
\end{equation}
Indeed, given $\sfT\in \GL_n(\R)$  define $\sfS\in \GL_n(\R)$ by $\sfS e_i=\|\sfT e_i\|_\X^{-1}\sfT e_i$  for   $i\in \n$. The triangle inequality for $\|\cdot\|_\X$ gives $\|\sfS\|_{\ell_1^n\to \X}\le 1$. As $\sfS^{-1} x=(\|\sfT e_1\|_\X\langle\sfT^{-1}x,e_1\rangle,\ldots,\|\sfT e_n\|_\X\langle \sfT^{-1}x,e_n\rangle)$ for every $x\in \R^n$, where $\langle\cdot,\cdot\rangle$ is the standard scalar product on $\R^n$, using Cauchy–Schwarz we get:
$$
\|\sfS^{-1} x\|_1=\sum_{i=1}^n \|\sfT e_i\|_\X|\langle \sfT^{-1}x,e_i\rangle|\le \bigg(\sum_{i=1}^n \|\sfT e_i\|_\X^2\bigg)^{\frac12}\|\sfT^{-1}x\|_2\le \bigg(\sum_{i=1}^n \|\sfT e_i\|_\X^2\bigg)^{\frac12}\|\sfT^{-1} \|_{\X\to \ell_2^n}\|x\|_\X.
$$ 
Thus, $\|\sfS^{-1}\|_{\X\to \ell_1^n}\le (\sum_{i=1}^n \|\sfT e_i\|_\X^2)^{1/2}\|\sfT^{-1} \|_{\X\to \ell_2^n}$.   Now~\eqref{eq:og bouds BM} follows by choosing $\sfT$ to be a minimizer in~\eqref{eq:def og}.

In light of~\eqref{eq:og bouds BM}, the rest of the discussion herein is devoted to bounding  $\og_\X$. For this, we will next rewrite $\og_\X$ as the solution of another minimization problem; one arrives at it unavoidably as it is designed so that its first variation coincides with the first variation of the   minimization  defining  $\og_\X$. This, of course, does not imply that the two minima must coincide, but it motivates the subsequent considerations. After the fact, we do not need to work herein with those first variations at all, as the ensuing proof relating the  two minimization problems  has a short justification using elementary linear algebra. 

For every $\sfT\in \sfM_n(\R)$ denote:
\begin{equation}\label{eq:def tau}
\tau_\X(\sfT)\eqdef \inf\bigg\{\Big(\sum_{i=1}^n s_i\Big)^{\frac12}:\ \exists s_1,\ldots,s_n\ge 0,\ \exists x_1,\ldots,x_n\in B_\X\ \mathrm{such\ that\ } \sfT\sfT^*=\sum_{i=1}^n s_ix_i\otimes x_i \bigg\}.
\end{equation}
This is well-defined as $\sfT\sfT^*=\sum_{i=1}^n (\sfT e_i)\otimes (\sfT e_i)$. Note that while~\eqref{eq:def tau} is reminiscent of the projective tensor product norm, restricting to $n$ summands in~\eqref{eq:def tau} makes $\tau_\X$  a substantially different object. In particular, $\tau_\X$ is positively homogeneous of order $1$ but it need not be a norm on $\sfM_n(\R)$. 

Consider two matrix sets $\Sigma_\X,\Lambda_\X\subset \sfM_n(\R)$ that are given as follows: 
\begin{equation}\label{eq:def Lambda}
\Sigma_\X\eqdef \big\{x\otimes x:\ x\in B_\X\big\}\qquad\mathrm{and}\qquad \Lambda_\X\eqdef \big\{\sfT\sfT^*:\ \sfT\in\sfM_n(\R)\ \mathrm{and}\ \tau_\X(\sfT)\le 1\big\}.
\end{equation}
With these notations, define the following parameter, which looks like the usual $K$-functional except that $\Lambda_\X$ need not be convex, and we are dealing with positive definite matrices: 
\begin{equation}\label{eq:def kappa}
\forall t\ge 0,\qquad \kappa_\X(t)=\min\Big\{
\|\sfA^{-\frac12}\|_{\X\to\ell_2^n}:
\sfA\in\big(\Lambda_\X+ t \conv\Sigma_\X\big)\cap\GL_n(\R)
\Big\}.
\end{equation}
The set in the right hand side of~\eqref{eq:def kappa} is nonempty as it contains a small enough positive multiple of the identity matrix $\Id_n\in \sfM_n(\R)$. The following lemma relates $\og_\X$ and $\kappa_\X$: 

\begin{lemma}\label{lem:relate the two optimizations}
 $\og_\X\sqrt{n}=\min\big\{\kappa_\X(t)+t\kappa_\X(t)^2: t\ge 0\big\}$.
\end{lemma}
It is mechanical to check that $\og_\X\sqrt{n}=\kappa_\X(0)$; we leave this as a straightforward exercise since it does not occur in the proof of Theorem~\ref{thm:main}. The direction of Lemma~\ref{lem:relate the two optimizations} that we will use is that $\og_\X\sqrt{n}\le \kappa_\X(t)+t\kappa_\X(t)^2$ for  $t>0$; as we will see in Section~\ref{sec:proof of lem} below, this is a fact in linear algebra. 

Passing to the derivation of Theorem~\ref{thm:main} from Lemma~\ref{lem:relate the two optimizations}, we now arrive at the only part of the reasoning herein  that is not algebraic manipulations, though it is natural in the local theory of Banach spaces, as it mimics a key step in the proofs  of e.g.~John's theorem~\cite{MR30135}, the Dvoretzky--Rogers lemma~\cite{DR50}, and the Lewis lemma~\cite{Lew79}. Fix $t>0$ and let $\sfA_t\in \sfM_n(\R)$  be a matrix that maximizes the determinant over the set appearing in the definition~\eqref{eq:def kappa} of $\kappa_\X(t)$, namely:
\begin{equation}\label{eq:choose A maximizer}
\det \sfA_t=\max \big\{
\det \sfA:
\sfA\in\Lambda_\X+ t \conv\Sigma_\X
\big\}.
\end{equation} 
Thus, while~\eqref{eq:def kappa} minimizes over invertible $\sfA\in \Lambda_\X+ t \conv\Sigma_\X$  the smallest $r>0$ for which $r\sfA^{1/2} B^n\supseteq B_\X$, we now consider such $\sfA$ that maximizes the volume of   $\sfA^{1/2} B^n$, which (following Auerbach and John) is well known to be a much better behaved relaxation of a containment constraint. 

As $\Lambda_\X+ t \conv\Sigma_\X$ has nonempty interior in  $\sfM_n(\R)$, we know that $\det \sfA_t>0$, whence $\sfA_t$ belongs to the set that appears in the right hand side of~\eqref{eq:def kappa}. By combining~\eqref{eq:og bouds BM} with the estimate $\og_\X\sqrt{n}\le \kappa_\X(t)+t\kappa_\X(t)^2$ of  Lemma~\ref{lem:relate the two optimizations}, and using the definition~\eqref{eq:def kappa} of $\kappa_\X(t)$, we therefore get:  
\begin{equation}\label{eq:use lemma}
d_{\BM}(\ell_1^n,\X)\le \big\|\sfA_t^{-\frac12}\big\|_{\X\to\ell_2^n} +t\big\|\sfA_t^{-\frac12}\big\|_{\X\to\ell_2^n}^2. 
\end{equation}

To bound the right hand side of~\eqref{eq:use lemma}, we will prove the following operator norm bound:  
\begin{equation}\label{eq:containment for optimizer}
\big\|\sfA_t^{-\frac12}\big\|_{\X\to\ell_2^n}  \le \frac{\sqrt{2n}}{\sqrt{\sqrt{4t+t^2}+t}}\le \frac{\sqrt{n}}{\sqrt[4]{t}}.
\end{equation}
A substitution of~\eqref{eq:containment for optimizer} into~\eqref{eq:use lemma} gives $d_{\BM}(\X,\ell_1^n)\le \sqrt{n}/\sqrt[4]{t}+n\sqrt{t}$,  yielding  Theorem~\ref{thm:main}  for $t\asymp n^{-2/3}$. The deduction of~\eqref{eq:containment for optimizer} appears in Section~\ref{sec:proof}. It consists of standard consequences of the fact that $\sfA_t$ is a  maximizer  per~\eqref{eq:choose A maximizer}, by examining how rank-one perturbations of $\sfA_t$ (corresponding to swaps of summands that occur in the representation of $\sfA_t$ as a member of $\Lambda_\X+ t \conv\Sigma_\X$) influence its determinant.

\section{Proof of Lemma~\ref{lem:relate the two optimizations}}\label{sec:proof of lem}

Fix $t\ge 0$. Here we will  prove that $\og_\X\sqrt{n}\le \kappa_\X(t)+t\kappa_\X(t)^2$. After an arbitrarily small perturbation to get invertible matrices, it suffices to show that if
 $\sfT,\sfC\in \GL_n(\R)$ are such that $\tau_\X(\sfT)\le 1$ and $\sfC\in \conv\Sigma_\X$, then: 
\begin{equation}\label{eq:hgamma goal}
\og_\X\sqrt{n}\le \big\|(\sfT\sfT^*+t\sfC)^{-\frac12}\big\|_{\X\to \ell_2^n} +t\big\|(\sfT\sfT^*+t\sfC)^{-\frac12}\big\|_{\X\to \ell_2^n}^2. 
\end{equation}
For~\eqref{eq:hgamma goal}, by the definition~\eqref{eq:def og} of $\og_\X$ it suffices to exhibit $\sfM\in \GL_n(\R)$ that satisfies: 
\begin{equation}\label{eq:two goals}
\Big(\sum_{i=1}^n \|\sfM e_i\|_\X^2\Big)^{\frac12}\le 1+t\big\|(\sfT\sfT^*+t\sfC)^{-\frac12}\big\|_{\X\to \ell_2^n}\qquad \mathrm{and}\qquad \|\sfM^{-1}\|_{\X\to \ell_2^n}=\big\|(\sfT\sfT^*+t\sfC)^{-\frac12}\big\|_{\X\to \ell_2^n}.
\end{equation}

As $\tau_\X(\sfT)\le 1$, there are $x_1,\ldots,x_n\in B_\X$ and $s_1,\ldots,s_n\ge 0$ with $\sum_{i=1}^n s_i\le 1$ and $\sfT\sfT^*=\sum_{i=1}^n s_i x_i\otimes x_i$. Define $\sfS\in \GL_n(\R)$ by $\sfS e_i=\sqrt{s_i} x_i$ for every $i\in \n$. Then:
\begin{equation}\label{eq:S identities}
\sfS\sfS^*=\sfT\sfT^*\qquad\mathrm{and}\qquad \sum_{i=1}^n \|\sfS e_i\|_\X^2\le 1.
\end{equation}  
From $\sfC\in \conv\Sigma_\X$ fix  $y_1,\ldots,y_m\in B_\X$ and  $\lambda_1,\ldots,\lambda_m>0$ with $\sum_{j=1}^m\lambda_j=1$ such that $\sfC=\sum_{j=1}^m \lambda_j y_j\otimes y_j$, i.e.,
\begin{equation}\label{eq:C identity}
\forall z\in \R^n,\qquad  z=\sum_{j=1}^m \lambda_j\langle y_j,\sfC^{-1} z\rangle y_j.
\end{equation}

In terms of the above data, the   $\sfM\in \GL_n(\R)$ for~\eqref{eq:two goals} is simply given by the following expression:
\begin{equation}\label{eq:our M}
\sfM\eqdef \sfS\big(\Id_n+t\sfS^{-1}\sfC(\sfS^*)^{-1}\big)^{\frac12}. 
\end{equation}
Then,
$$
\sfM\sfM^*\stackrel{\eqref{eq:our M}}{=}\sfS\big(\Id_n+t\sfS^{-1}\sfC(\sfS^*)^{-1}\big)\sfS^*=\sfS\sfS^*+t\sfC\stackrel{\eqref{eq:S identities}}{=}\sfT\sfT^*+t\sfC. 
$$
So, $\|\sfM^{-1} x\|_2=\sqrt{\langle (\sfM\sfM^*)^{-1}x,x\rangle} =\sqrt{\langle (\sfT\sfT^*+t\sfC)^{-1}x,x\rangle} =\|(\sfT\sfT^*+t\sfC)^{-\frac12}x\|_2$ for every $x\in \R^n$, whence the second part of~\eqref{eq:two goals} indeed holds.

We will next establish the following estimate: 
\begin{equation}\label{eq:two goals'}
\Big(\sum_{i=1}^n \|(\sfM-\sfS) e_i\|_\X^2\Big)^{\frac12}\le t\big\|(\sfT\sfT^*+t\sfC)^{-\frac12}\big\|_{\X\to \ell_2^n}.
\end{equation}
Note that \eqref{eq:two goals'}  implies the first part of~\eqref{eq:two goals} by the triangle inequality and the second part of~\eqref{eq:S identities}. Because 
$$
\sfS^{-1}(\sfM-\sfS)\Big(\big(\Id_n+t\sfS^{-1}\sfC(\sfS^*)^{-1}\big)^{\frac12}+\Id_n\Big)\stackrel{\eqref{eq:our M}}{=}\Big(\big(\Id_n+t\sfS^{-1}\sfC(\sfS^*)^{-1}\big)^{\frac12}-\Id_n
\Big)\Big(\big(\Id_n+t\sfS^{-1}\sfC(\sfS^*)^{-1}\big)^{\frac12}+\Id_n\Big)=t\sfS^{-1}\sfC(\sfS^*)^{-1},
$$
we have:
$$
\sfC^{-1}(\sfM-\sfS)= t(\sfS^*)^{-1} \Big(\Id_n+\big(\Id_n+t\sfS^{-1}\sfC(\sfS^*)^{-1}\big)^{\frac12}\Big)^{-1}.
$$
Consequently, letting $\preceq$ denote the positive semidefinite ordering, we have:
\begin{multline*}
\big(\sfC^{-1}(\sfM-\sfS)\big)\big(\sfC^{-1}(\sfM-\sfS)\big)^*=t^2(\sfS^*)^{-1} 
\Big(\Id_n+\big(\Id_n+t\sfS^{-1}\sfC(\sfS^*)^{-1}\big)^{\frac12}\Big)^{-2}\sfS^{-1}\\ \preceq t^2(\sfS^*)^{-1} 
\big(\Id_n+t\sfS^{-1}\sfC(\sfS^*)^{-1}\big)^{-1}\sfS^{-1}=t^2(\sfS\sfS^*+t\sfC)^{-1}\stackrel{\eqref{eq:S identities}}{=}t^2(\sfT\sfT^*+t\sfC)^{-1},
\end{multline*}
which is equivalent to the following statement: 
\begin{equation}\label{eq:PSD order}
\forall x\in \R^n,\qquad \|(\sfM-\sfS)^*\sfC^{-1}x\|_2\le t\big\|(\sfT\sfT^*+t\sfC)^{-\frac12}x\big\|_2\le t\big\|(\sfT\sfT^*+t\sfC)^{-\frac12}\big\|_{\X\to \ell_2^n}\|x\|_\X. 
\end{equation}
From here we conclude the justification of~\eqref{eq:two goals'} as follows: 
\begin{align*}
\Big(\sum_{i=1}^n &\|(\sfM-\sfS) e_i\|_\X^2\Big)^{\frac12}\stackrel{\eqref{eq:C identity}}{=} \Big(\sum_{i=1}^n \|\sum_{j=1}^m\lambda_j\langle (\sfM-\sfS)^*\sfC^{-1} y_j,e_i\rangle y_j\|_\X^2\Big)^{\frac12}
\le \sum_{j=1}^m\lambda_j \|y_j\|_\X\Big(\sum_{i=1}^n \langle  (\sfM-\sfS)^*\sfC^{-1} y_j, e_i\rangle^2 \Big)^{\frac12}\\&=\sum_{j=1}^m\lambda_j \|y_j\|_\X \|(\sfM-\sfS)^*\sfC^{-1} y_j\|_2 \stackrel{\eqref{eq:PSD order}}{\le}  t\big\|(\sfT\sfT^*+t\sfC)^{-\frac12}\big\|_{\X\to \ell_2^n}\sum_{j=1}^m\lambda_j \|y_j\|_\X^2\le t\big\|(\sfT\sfT^*+t\sfC)^{-\frac12}\big\|_{\X\to \ell_2^n}. \tag*{\qed}
\end{align*}

\section{Deduction of Theorem~\ref{thm:main} from Lemma~\ref{lem:relate the two optimizations}}\label{sec:proof}

Fix $t>0$. Recall that in the Introduction we defined $\sfA_t\in \GL_n(\R)$ by~\eqref{eq:choose A maximizer}, and saw that, using Lemma~\ref{lem:relate the two optimizations}, for Theorem~\ref{thm:main} it suffices to prove that~\eqref{eq:containment for optimizer} holds. For this,  as $\sfA_t$ belongs to the compact set $\Lambda_\X+ t \conv\Sigma_\X$, recalling~\eqref{eq:def tau} and~\eqref{eq:def Lambda} we can write  $\sfA_t=\sum_{i=1}^n s_i x_i\otimes x_i+ t\sum_{k=1}^m \lambda_ky_k\otimes y_k$ for some   $x_1,\ldots,x_n,y_1,\ldots,y_m\in B_\X$ as well as $s_1,\ldots,s_n\ge 0$ with $\sum_{i=1}^ns_i\le 1$ and $\lambda_1,\ldots,\lambda_m>0$ with $\sum_{k=1}^m\lambda_k=1$.  Without loss of generality, $s_1,\ldots,s_n>0$ and $\sum_{i=1}^ns_i=1$. Indeed, we cannot have $s_1=\ldots=s_n=0$ as otherwise for any $x\in B_\X\setminus \{0\}$ the matrix $x\otimes x+\sfA_t$ belongs to $\Lambda_\X+ t \conv\Sigma_\X$ and its determinant is larger than $\det \sfA_t$. So, $S=\sum_{i=1}^n s_i>0$. We cannot have  $\sum_{i=1}^n s_i< 1$ as then $\sfA_t+(1/S-1)\sum_{i=1}^n s_ix_i\otimes x_i\in \Lambda_\X+ t \conv\Sigma_\X$  and its determinant is larger than $\det\sfA_t$. If $\ell=|\{i\in \n:\ s_i>0\}|<n$, then for any $i\in \n$ for which $s_i>0$ we can replace the summand $s_i x_i\otimes x_i$ by $n-\ell+1$ copies of $(s_i/(n-\ell+1))x_i\otimes x_i$ to get $s_1,\ldots,s_n> 0$.

Fix $j\in \{1,\ldots,m\}$ and $z\in B_\X$. As $\sfA_t+th(z\otimes z-y_j\otimes y_j)\in \Lambda_\X+ t \conv\Sigma_\X$ for every $0\le h\le \lambda_j$, we have:
\begin{align*}
0& \stackrel{\eqref{eq:choose A maximizer}}{\ge}  \left.\frac{\ud }{\ud h}\right|_{h=0^+} \det\big(\sfA_t+th(z\otimes z-y_j\otimes y_j)\big)
=t \Big(\langle \sfA_t^{-1} z,z\rangle-\langle \sfA_t^{-1} y_j,y_j\rangle\Big)\det\sfA_t.
\end{align*}
Taking the supremum over $z\in B_\X$, we obtain: 
\begin{equation}\label{eq:equality on the yj} 
\forall j\in \{1,\ldots,m\},\qquad \big\|\sfA_t^{-\frac12} y_j\big\|_2=\big\|\sfA_t^{-\frac12} \big\|_{\X\to \ell_2^n}. 
\end{equation}
Also, for  distinct $i,k\in \n$, if $|h|< \min\{s_i,s_k\}$, then  $\sfA_t+h(x_i\otimes x_i-x_k\otimes x_k)\in \Lambda_\X+ t \conv\Sigma_\X$, so we get:
$$
0 \stackrel{\eqref{eq:choose A maximizer}}{=}  \left.\frac{\ud }{\ud h}\right|_{h=0} \det\big(\sfA_t+h(x_i\otimes x_i-x_k\otimes x_k)\big)=\Big(\big\|\sfA_t^{-\frac12} x_i\big\|_2^2-\big\|\sfA_t^{-\frac12} x_k\big\|_2^2\Big)\det\sfA_t.
$$
Consequently, there exists $\alpha>0$ such that: 
\begin{equation}\label{eq:equality on the xj} 
\forall i\in \n,\qquad \big\|\sfA_t^{-\frac12} x_i\big\|_2=\alpha. 
\end{equation}
Using \eqref{eq:equality on the yj} and \eqref{eq:equality on the xj}, we obtain the following identity: 
$$
n=\Tr\big(\sfA_t^{-1}\sfA_t\big)=\sum_{i=1}^n s_i\Tr\big(\sfA_t^{-1}x_i\otimes x_i\big)+t\sum_{j=1}^m \lambda_j \Tr\big(\sfA_t^{-1}y_j\otimes y_j\big)
= \alpha^2 +t\big\|\sfA_t^{-\frac12} \big\|_{\X\to \ell_2^n}^2,
$$
i.e.,
\begin{equation}\label{eq:n trace formula like john}
\alpha^2=n-t\big\|\sfA_t^{-\frac12} \big\|_{\X\to \ell_2^n}^2.
\end{equation}

Fix $i_0\in \n$ with $s_{i_0}\le 1/n$, so $s_{i_0}\alpha^2<1$ by~\eqref{eq:n trace formula like john}. Whence, $\sfA_t-s_{i_0}x_{i_0}\otimes x_{i_0}\succeq0$ is invertible because:
$$
\det(\sfA_t-s_{i_0}x_{i_0}\otimes x_{i_0})=\Big(1-s_{i_0}\Tr\big(\sfA_t^{-1}x_{i_0}\otimes x_{i_0}\big)\Big)\det\sfA_t=\Big(1-s_{i_0}\big\|\sfA_t^{-\frac12} x_{i_0}\big\|_2^2\Big)\det\sfA_t\stackrel{\eqref{eq:equality on the xj} }{=}(1-s_{i_0}\alpha^2)\det\sfA_t\neq 0.
$$
Since  $\sfA_t-s_{i_0}x_{i_0}\otimes x_{i_0}+s_{i_0} x\otimes x\in \Lambda_\X+ t \conv\Sigma_\X$ for every $x\in B_\X$, we therefore get: 
\begin{align*}
 \Big(1+s_{i_0}\big\langle &(\sfA_t-s_{i_0}x_{i_0}\otimes x_{i_0})^{-1}x_{i_0},x_{i_0}\big\rangle\Big)\det\big(\sfA_t-s_{i_0}x_{i_0}\otimes x_{i_0}\big)=\det \sfA_t\\& \stackrel{\eqref{eq:choose A maximizer}}{\ge} \det\big(\sfA_t-s_{i_0}x_{i_0}\otimes x_{i_0}+s_{i_0} x\otimes x\big)=\Big(1+s_{i_0}\big\langle (\sfA_t-s_{i_0}x_{i_0}\otimes x_{i_0})^{-1}x,x\big\rangle\Big)\det\big(\sfA_t-s_{i_0}x_{i_0}\otimes x_{i_0}\big).
\end{align*}
Hence, we have:
\begin{align}\label{eq:norm attained at x_0}
\begin{split}
\big\|(\sfA_t-s_{i_0}x_{i_0}\otimes x_{i_0})^{-\frac12}x_{i_0}\big\|_2^2&=\big\langle (\sfA_t-s_{i_0}x_{i_0}\otimes x_{i_0})^{-1}x_{i_0},x_{i_0}\big\rangle\\&\ge \big\langle (\sfA_t-s_{i_0}x_{i_0}\otimes x_{i_0})^{-1}x,x\big\rangle=\big\|(\sfA_t-s_{i_0}x_{i_0}\otimes x_{i_0})^{-\frac12}x\big\|_2^2.
\end{split}
\end{align}
Because this holds for every $x\in B_\X$ and $\sfA_t^{-1}\preceq (\sfA_t-s_{i_0}x_{i_0}\otimes x_{i_0})^{-1}$, we conclude that
\begin{align}\label{eq:before rank one inverse}
\begin{split}
\|\sfA_t^{-\frac12}\big\|_{\X\to\ell_2^n}^2&\le \big\|(\sfA_t-s_{i_0}x_{i_0}\otimes x_{i_0})^{-\frac12}\big\|_{\X\to \ell_2^n}^2\\& \stackrel{\eqref{eq:norm attained at x_0}}{=}\big\|(\sfA_t-s_{i_0}x_{i_0}\otimes x_{i_0})^{-\frac12}x_{i_0}\big\|_2^2=\big\langle (\sfA_t-s_{i_0}x_{i_0}\otimes x_{i_0})^{-1}x_{i_0},x_{i_0}\big\rangle.
\end{split}
\end{align}

Finally, the Sherman--Morrison rank-one inverse formula (see e.g.~\cite[equation~(0.7.4.2)]{HJ13}) gives: 
\begin{equation}\label{eq:apply rank one inverse}
\big\langle
(\sfA_t-s_{i_0} x_{i_0}\otimes x_{i_0})^{-1}x_{i_0},x_{i_0}
\big\rangle
=
\frac{\langle\sfA_t^{-1}x_{i_0},x_{i_0}\rangle}
{1-s_{i_0}\langle\sfA_t^{-1}x_{i_0},x_{i_0}\rangle}
\stackrel{\eqref{eq:equality on the xj} }{=}\frac{\alpha^2}{1-s_{i_0}\alpha^2}
\leq\frac{\alpha^2}{1-\alpha^2/n}
\stackrel{\eqref{eq:n trace formula like john}}{=}
\frac{n^2}
{t\big\|\sfA_t^{-\frac12}\big\|_{\X\to\ell_2^n}^2}-n.
\end{equation}
A substitution of~\eqref{eq:apply rank one inverse} into~\eqref{eq:before rank one inverse} and solving the resulting quadratic inequality gives~\eqref{eq:containment for optimizer}. \qed

\begin{remark}\label{rem:gluskin} {\em Fix $1\le \alpha\le \sqrt{n}$. In the Introduction we stated that~\cite{Glu81} provides a random normed space $\mathbf{G}_\alpha=(\R^n,\|\cdot\|_{\mathbf{G}_\alpha})$ such that $\sR_{\BM}(\mathbf{G}_\alpha)\asymp \alpha\sqrt{n}$ with high probability. What~\cite{Glu81} actually gets is a random normed space $\mathbf{G}=(\R^n,\|\cdot\|_{\mathbf{G}})$ such that if $\mathbf{G'}$ is an independent copy of $\mathbf{G}$, then  with high probability $d_{\BM}(\mathbf{G'},\mathbf{G})\asymp n$. Let $\mathcal{E}\subset \R^n$ be a (random) ellipsoid such that $\mathcal{E}\subset B_{\mathbf{G}}\subset d_{\BM}(\ell_2^n,\mathbf{G})\mathcal{E}$. Define $\mathbf{G}_\alpha$ to be the normed space whose unit ball is $B_{\mathbf{G}} \cap (\alpha\mathcal{E})$.  As $\alpha\ge 1$, we have $\mathcal{E}\subset B_{\mathbf{G}} \cap(\alpha \mathcal{E})\subset \alpha\mathcal{E}$, so $d_{\BM}(\ell_2^n,\mathbf{G}_\alpha)\le \alpha$. Also,  $\min\{\alpha/d_\BM(\ell_2^n,\mathbf{G}),1\}B_{\mathbf{G}}\subset B_{\mathbf{G}} \cap(\alpha \mathcal{E})\subset B_{\mathbf{G}}$, so $d_{\BM}(\mathbf{G}_\alpha,\mathbf{G})\le \max\{d_\BM(\ell_2^n,\mathbf{G})/\alpha,1\}\le \sqrt{n}/\alpha$, where we used $d_\BM(\ell_2^n,\mathbf{G})\le \sqrt{n}$ (by John's theorem) and  $\alpha\le \sqrt{n}$. Consequently, with high probability we have:
$$
\alpha\sqrt{n}\ge d_{\BM}(\ell_2^n,\mathbf{G}_\alpha) \sR_{\BM}(\ell_2^n)\ge \sR_{\BM}(\mathbf{G}_\alpha)\ge  d_{\BM}(\mathbf{G}',\mathbf{G}_\alpha)\ge\frac{d_{\BM}(\mathbf{G'},\mathbf{G})}{d_{\BM}(\mathbf{G}_\alpha,\mathbf{G})}\gtrsim \frac{n}{\sqrt{n}/\alpha}=\alpha\sqrt{n}. 
$$
 }
\end{remark}

\bibliographystyle{alphaabbrvprelim}
\bibliography{gamma}

\end{document}